\documentclass[12pt]{article}
\def\capa{{\mathrm{cap}\,}}
\def\diam{{\mathrm{diam}\,}}
\def\C{{\mathbf{C}}}
\def\id{{\mathrm{id}}}
\def\const{{\mathrm{const}}}
\author{Alexandre Eremenko\thanks{Supported by NSF grants
DMS-0100512 and DMS-0244421.}}
\title{A Markov-type inequality for arbitrary plane continua}
\begin{document}
\maketitle
\begin{abstract}
Markov's inequality is
$$\sup_{[-1,1]}|f'|\leq(\deg f)^2\sup_{[-1,1]}|f|,$$
for all polynomials $f$.
We prove a precise version of this inequality with
an arbitrary continuum in the complex plane instead
of the interval $[-1,1]$.
\end{abstract}

\noindent
{\bf Theorem 1.} {\em Let $E$ be a continuum
in the complex plane,
and $f$ a polynomial of degree $d$. Then
\begin{equation}\label{1}
\capa E\sup_E|f'|\leq 2^{1/d-1}d^2\sup_E|f|.
\end{equation}
}

Here $\capa$ denotes the transfinite diameter (capacity)
of a set \cite[Ch. 2]{Ahlfors}. Using the well-known
inequality $\diam E\leq 4\,\capa E$ we obtain
\vspace{.1in}

\noindent
{\bf Corollary.} {\em With the assumptions
of Theorem 1 we have
\begin{equation}
\label{6}
\diam E\sup_E|f'|\leq 2^{1/d+1}d^2\sup_E|f|.
\end{equation}
}
This inequality looks more elementary than (\ref{1}) but we could not
find a direct proof.

To compare these results with Markov's inequality \cite{Mark},
\cite{Markoff}, we take $E=[-1,1]$.
Then $\capa E=1/2$, and $\diam E=2$, so (\ref{1}) and
(\ref{6})
become 
$$\sup_{[-1,1]}|f'|\leq 2^{1/d} d^2\sup_{[-1,1]}|f|,$$
while Markov's inequality is
$$\sup_{[-1,1]}|f'|\leq d^2\sup_{[-1,1]}|f|.$$
Thus the estimate of Theorem 1 for $E=[-1,1]$
is asymptotically
exact as $d\to\infty$.

Both (\ref{1}) and Markov's inequality
are precise, and in both of them equality
is achieved on the Chebyshev polynomials,
but to obtain equality in (\ref{1}) one has to
take $E$ to be the level set of the Chebyshev polynomial,
$E=\{ z:|T_d(z)|\leq 1\}$, rather than the interval $[-1,1]$.

Inequality (\ref{6}) for $E=[-1,1]$ is
worse than Markov's inequality by the factor
$2^{1/d}$, and one can show that equality never
takes place in (\ref{6}).

Pommerenke \cite{Po} proved a weaker form of (\ref{1}):
$$\capa E\sup_E|f'|\leq (ed^2/2)\sup_E|f|.$$
A special case of (\ref{1}) is the following result
conjectured by Erd\H{o}s and proved in \cite{EL}:
\vspace{.1in}

\noindent
{\bf Theorem A.} {\em Let $f$ be a monic polynomial such that the level set
$E(f)=\{ z:|f(z)|\leq 1\}$ is connected. Then 
$$\sup_{E(f)}|f'|\leq 2^{1/d-1}d^2.$$
}
\vspace{.1in}

Equality in Theorem A occurs if and only if 
$f(z)=c^{-d}T_d(2^{1/d-1}cz+b)$, where $|c|=1$, and
$T_d$ is the Chebyshev polynomial defined by
$$\cos dz=T_d(\cos z).$$

Theorem A will be used in the proof of the general
result. Theorem~1 evidently follows from 
\vspace{.1in}

\noindent
{\bf Theorem 2.} {\em Let $f$ be a polynomial of degree $d$,
and $E$
a component of the set $E(f)=\{ z:|f(z)|\leq 1\}$. 
Then
\begin{equation}
\label{8}
\capa E\sup_E|f'|\leq 2^{1/d-1}d^2.
\end{equation}
}
\vspace{.1in}

This estimate is best possible, and the equality holds
if and only if $f(z)=aT_d(cz+b)$ with $c\neq 0,\, |a|=1,$ and
$E=E(f)$.

The proof of Theorem 2~is based on the idea which
was used in
\cite{ELev,EH} to prove other polynomial
inequalities.
The
method\footnote{
A method is a device which
you used twice \cite{Polya}.} consists in
including $f$ in a holomorphic
family of polynomials parametrized by critical values,
and proving that the quantities of interest are subharmonic
functions of parameters. Then it remains to check the inequality
on the boundary of the parameter domain.

Following \cite{EH}, we begin by establishing existence
of extremal polynomials. Then we will show that for every
monic extremal polynomial, the set $E(f)$ is connected,
thus reducing Theorem~2 to Theorem~A.
\vspace{.1in}

\noindent
{\bf Lemma 1.} {\em For every positive integer
$d$ there exists a
polynomial $f^*$ of degree at most $d$
and a component $E^*$ of the set $E(f^*)$
such that
\begin{equation}
\label{2}
\capa E^*\sup_{E^*}|(f^*)^\prime|\geq \capa E\sup_E|f'|,
\end{equation}
for every polynomial $f$ of degree at most $d$ and every component
$E$ of $E(f)$.}
\vspace{.1in}

{\em Proof}. Let $(f_{n})$ and $(E_n)$, be extremizing
sequences.
Notice that the product $\capa E_n\sup_{E_n}|f^\prime_n|$
does not change if the independent variable $z$ is replaced by $az+b$,
where $a\neq 0$. So we can use the normalization
\begin{equation}\label{3}
\capa E_n=1,\quad\mbox{and} \quad 0\in E_n.
\end{equation}
This implies $\capa E(f_n)\geq 1$, so
\begin{equation}\label{17}
f_n(z)=a_nz^d_n+\ldots,\quad\mbox{where}\quad |a_n|\leq 1,\quad\mbox{and}
\quad d_n\leq d.
\end{equation}
Let $g_n$ be the Green function
of the complement of $E_n$ with
the pole at infinity. We have the Riesz representation
$$g_n(z)=\int_\C\log|z-\zeta|\,d\mu_n(\zeta),$$
with some probability measures $\mu_n$ whose supports belong to
$E_n$. Any continuum of capacity $1$ has diameter at most $4$,
\cite{Ahlfors},
so our normalization conditions (\ref{3})
imply that the supports of $\mu_n$ are contained in the
disc
$D(4)=\{ z:|z|\leq 4\}$. 
Thus one can choose a subsequence
such that $\mu_n\to\mu$, where $\mu$
is a probability measure
on $D(4)$. Then we have $g_n\to g$ uniformly in the region
$\Delta(5)=\{ z:|z|>5\}$,
so $g_n(z)\leq 2g(z),\; z\in \Delta(5)$,
and by the Principle of
Harmonic Majorant, $|f_n(z)|<\exp(2dg(z)),\;
z\in \Delta(5)$, where we used $|a_n|\leq 1$ and $d_n\leq d$ from
(\ref{17}). Thus
the sequence of our polynomials
is uniformly bounded on compact subsets of
the plane and so, after choosing a subsequence, 
we obtain a limit polynomial $f^*$.
Now (\ref{2}) is evident.
This completes the proof of the Lemma.
\vspace{.1in}

We continue the proof of Theorem 2.
For fixed $d$, we will call $(f^*,E^*)$ an extremal pair,
if $\deg f^*\leq d$, $E^*$ is a component of $E(f^*)$
and (\ref{2}) is satisfied.

To simplify our notation, let $(f,E)$
be an extremal pair.
If the set $E(f)$ is connected that is $E=E(f)$,
then we can normalize as in 
(\ref{3}), and the polynomial $f$ will be monic
because $\capa E(f)=1$.
So Theorem~A applies and (\ref{8}) follows.
Now let us assume that the set $E(f)$ is disconnected.

This assumption implies that there is a critical value $a$ of $f$
such that $|a|>1$. Let $D$ be a disc centered at the point $a$
of radius $8\delta<|a|-1$, so that it does not intersect
the unit disc. 
Let $\chi$ be a continuously differentiable function
in the complex plane $\C$, with support in $D$,
such that $\chi(a)=1$ and
$|\nabla\chi|<1/(2\delta)$. An example
of such function is $(1-|z-a|^2/(64\delta^2))^2$.
Let $\Delta$ be the disc of radius $\delta$ centered at $0$.
Then for 
$\lambda\in\Delta,$ the map
$$
\psi_\lambda=\id+\lambda.\chi$$
is a quasiconformal homeomorphism of the plane, holomorphic
outside $D$ and holomorphically depending on $\lambda$.
Now we construct a quasiconformal homeomorphism $\phi_\lambda$
which makes the function
\begin{equation}
\label{flambda}
f_\lambda=\psi_\lambda\circ f\circ\phi_\lambda^{-1}
\end{equation}
analytic in the whole plane. To achieve this, $\phi_\lambda$
has to satisfy a Beltrami equation (see \cite[Ch. 1C]{Ahlforsq}):
$$\frac{\partial\phi_\lambda}{\partial\overline{z}}=\mu_\lambda
\frac{\partial\phi_\lambda}{\partial z},$$
where the Beltrami coefficient 
$$\mu_\lambda=(\mu_{\psi_\lambda}\circ f)\frac{\overline{f'}}{f'},
\quad\mbox{and}\quad \mu_{\psi_\lambda}=\frac{\partial\psi_\lambda}{\partial\overline{w}}:\frac{\partial\psi_\lambda}{\partial w}.$$
Differentiation of (\ref{flambda}) with respect to
$\overline{\lambda}$ and application of the Beltrami equation shows
that $f_\lambda$ is analytic.
Notice that $\mu_{\psi_\lambda}$ and thus $\mu_\lambda$ depend
analytically on $\lambda$. The existence theorem for
Beltrami equation (also known as the ``Measurable Riemann Theorem'',
\cite[Ch. 2, thms. 4, 5]{Ahlfors})
says that there exists a unique homeomorphic solution
$\phi_\lambda:\C\to\C$
normalized by $\phi_\lambda(0)=0$ and $\phi_\lambda(1)=1$.
Moreover, this normalized solution depends analytically on $\lambda$.

As $f_\lambda$ defined in (\ref{flambda}) is analytic in $\C$, and
both $\psi_\lambda$ and $\phi_\lambda$ are homeomorphisms, we conclude
that $f_\lambda$ are polynomials of the same degree as $f$.

Now we claim that
$f_\lambda$ depends holomorphically on $\lambda$. More precisely,
\begin{equation}\label{family}
(z,\lambda)\mapsto f_\lambda(z), \quad\C\times\Delta\to\C
\end{equation}
is a holomorphic function of two variables.
To verify this, we apply the operator
$\partial/\partial\overline{\lambda}$ to both
sides of the equation
$$f_\lambda\circ\phi_\lambda=\psi_\lambda\circ f.$$
As $\psi_\lambda$ is holomorphic in $\lambda$, we obtain
$$\frac{\partial f_{\lambda}}{\partial\overline{\lambda}}\circ\phi_\lambda=
-\left(\frac{\partial f_\lambda}{\partial z}\circ\phi_\lambda\right)\frac{\partial\phi_\lambda}{\partial\overline{\lambda}}=0;$$
the last equality holds
because $\phi_\lambda$ depends holomorphically on $\lambda$. This proves that the map
(\ref{family}) is holomorphic.

Notice that the functions $z\mapsto\phi_\lambda(z)$ are holomorphic on $E$.
This is because the $\psi_\lambda$ are holomorphic in $\C\backslash D$
and the full preimage $f^{-1}(D)$
is disjoint from $E$ by our choice of $D$.

The polynomials $f_\lambda$ are obtained from $f$ by moving the critical
value $a$ of $f$ to the point $a+\lambda$. 

We put $E_\lambda=\phi_\lambda(E)$. Then $E_\lambda$ is a component
of the set $E(f_\lambda)$. Such families
of sets $E_\lambda$ parametrized by a complex parameter
$\lambda$
are called holomorphic motions; they were much studied recently
(see, for example, \cite{Holomo}), but
we don't use any deep properties of holomorphic motions here.
The following two lemmas are almost evident:
\vspace{.1in}

\noindent
{\bf Lemma 2.} {\em The function
$\lambda\mapsto \log\sup_{E_\lambda}|f^\prime_\lambda|$
is subharmonic.}
\vspace{.1in}

{\em Proof.} 
In some neighborhood $V$ of $E$ we have
$$f_\lambda\circ\phi_\lambda=f.$$
Differentiation with respect to $z$ for $z\in V$ 
gives
\begin{equation}
\label{hh}
(f^\prime_\lambda\circ\phi_\lambda(z))
\phi^\prime_\lambda(z)=
f'(z),\quad z\in E,
\end{equation}
where we used the fact that $\phi_\lambda$ is holomorphic
on $E$. It follows from (\ref{hh}) that the functions
$\lambda\mapsto f^\prime_\lambda\circ\phi_\lambda(z)$
are holomorphic for every fixed $z\in V$. So
$f^\prime_\lambda\circ\phi_\lambda$ is a holomorphic
function of two variables, $\lambda$ and $z$, for $z$ in
$V$ and sufficiently small $\lambda$. Taking logarithm
of the modulus and then taking supremum with respect
to $z\in E$ we arrive at the conclusion of the Lemma.
\vspace{.1in}

\noindent
{\bf Lemma 3.} {\em The function
$\lambda\mapsto \log\capa E_\lambda$ is subharmonic.}
\vspace{.1in}

This should be well-known and the proof below applies to
arbitrary holomorphic motion of a compact set,
but we could not find a convenient reference.

{\em Proof.}
Let
$$d_n(\lambda)=
\sup\prod_{1\leq i<j\leq n}|\phi_\lambda(z_i)-
\phi_\lambda(z_j)|^{2/n(n-1)},$$
where the supremum is taken over all subsets $\{ z_j\}$
of $E$ of cardinality $n$. Then $\log d_n(\lambda)$
is a subharmonic function. Now, $(d_n)$ is a
decreasing sequence \cite{Ahlfors}, and the
transfinite diameter can be defined by the formula
$$\capa E_\lambda=\lim_{n\to\infty}d_n,$$
which completes the proof.
\vspace{.1in}

So the function
$$F(\lambda)=
\log(\capa E_\lambda\sup_{E_\lambda}|f^\prime_\lambda|)$$ 
is subharmonic. As the pair $(f,E)$ is assumed to be
extremal, we obtain from the Maximum Principle
that 
\begin{equation}\label{const}
F=\const.
\end{equation}
To summarize, we proved that if $(f,E)$ is an extremal pair,
and $f$ has a critical value $a$ outside of the unit disc,
then for every $\lambda\in\Delta$,
all pairs $(f_\lambda,E_\lambda)$ are also extremal.
Notice that $f_\lambda$ has a critical value $a+\lambda$,
so taking $\lambda=\delta a/(2|a|)$ we obtain an extremal pair
where the polynomial has critical value $a(1+\delta/(2|a|))$.

By iterating this construction, we obtain a sequence of
extremal pairs $(f_k,E_k)$ where for each $k$ the polynomial $f_k$
has a critical value of modulus $|a|+k\delta/2$. 
  
We can normalize each $f_k$ and $E_k$
as in (\ref{3}),
so that the family $\{ f_k\}$ will be uniformly
bounded on compact subsets of $\C$,
and we will obtain a limit
pair $(f_\infty,E_\infty)$, where $f_\infty$ evidently
has smaller degree than $f$. From (\ref{const}) we conclude
that
\begin{equation}\label{hu}
\capa E\sup_E|f'|=\capa E_\infty\sup_{E_\infty}|f'_\infty|.
\end{equation}
Now we use induction on $d$.
Denote the left hand side of (\ref{2}) by $C(d)$.
Evidently $C(1)=1$. So Theorem 2 holds with $d=1$.
Suppose that $k>1$ is the smallest integer
for which the inequality
of Theorem 2 does not hold with $d=k$. Applying Lemma 1 with $d=k$ we obtain
an extremal pair $(f,E)$, where $\deg f=k$, because the right hand
side of (\ref{8}) is strictly increasing. Then 
\begin{equation}
\label{7}
C(k)>C(k-1),
\end{equation}
again because the right hand side of (\ref{8}) is strictly increasing. 
If $E(f)$ is connected then we use Theorem~A
to conclude that Theorem 2 holds with $d=k$,
contrary to our assumption.

If $E(f)$ is disconnected
we obtain a contradiction between (\ref{hu}) and
(\ref{7}), because $\deg f_\infty<k$ and Theorem 2
holds for all $d<k$. This contradiction shows that the number $k$
does not exist, so Theorem 2 holds for all $d$.
\vspace{.1in}

It remains to address the case of equality.
We already saw in the course of the proof that
equality cannot occur when $E(f)$ is disconnected.
This reduces the case of equality in Theorem 2 to 
the case of equality in Theorem~A. 

This paper was written when the author
was visiting 
Christian-Albrechts University of Kiel. 
The author thanks Vladimir Andrievski and Walter Bergweiler
for stimulating
discussions, and the referee for suggestions that
improved the exposition.

{\em Purdue University, West Lafayette, IN 47907 USA

eremenko@math.purdue.edu}
\end{document}